\documentclass[12pt]{article}

\usepackage{fullpage,amsfonts,amsmath,amsthm,amssymb,mathrsfs,graphicx,upgreek}
\usepackage{verbatim,url}
\usepackage{graphicx,hyperref,enumitem}
\usepackage[noadjust,sort]{cite}
\usepackage{tocloft}

   \usepackage[top=2.54cm, bottom=2.54cm, left=2.54 cm, right=2.54 cm]{geometry}
 \newcommand{\lab}[1]{\label{#1}}                

 \usepackage[usenames,dvipsnames]{color}

\newcommand{\jt}[1]{{\color{black}#1}}
\newcommand{\jtb}[1]{{\color{black}#1}}

\newcommand{\remove}[1]{}
\newcommand\eqn[1]{(\ref{#1})}

\newcommand{\be}{\begin{equation}}
\newcommand{\bel}[1]{\begin{equation}\lab{#1}\ }
\newcommand{\ee}{\end{equation}}
\newcommand{\bea}{\begin{eqnarray}}
\newcommand{\eea}{\end{eqnarray}}
\newcommand{\bean}{\begin{eqnarray*}}
\newcommand{\eean}{\end{eqnarray*}}

\newtheorem{thm}{Theorem}
\newtheorem{cor}[thm]{Corollary}
\newtheorem{con}[thm]{Conjecture}
\newtheorem{lemma}[thm]{Lemma}
\newtheorem{definition}[thm]{Definition}
\newtheorem{claim}[thm]{Claim}

\def\proof{\noindent{\em Proof.\ } }
\def\qed{~~\vrule height8pt width4pt depth0pt}

\newcommand{\ind}[1]{1_{\{#1\}}}

\def\ham{{\tt HAM}}
\def\aham{{\tt HAM^+}}
\def\connect{{\tt CNT}}
\def\length{{\tt LC}}
\def\expand{{\tt EXPN}}
\def\expandd{{\tt SEXPN}}
\def\expanddd{{\tt SSEXPN}}
\def\col{{\tt COL}}
\def\degree{{\tt D2}}
\def\typical{{\tt TPCL}}

\def\G{{\mathcal G}}

\def\calN{{\mathcal N}}

\def\calP{{\mathcal P}}

\def\calX{{\mathcal X}}
\def\calY{{\mathcal Y}}
\def\calZ{{\mathcal Z}}

\def\barp{{\bar p}}

\def\sP{{\mathscr P}}
\def\sE{{\mathscr E}}

\def\ex{{\mathbb E}}
\def\pr{{\mathbb P}}

\def\Bin{{\bf Bin}}

\def\bfn{{\bf n}}

\def\eps{\epsilon}

\def\ss{\smallskip}
\def\non{\nonumber}
\def\no{\noindent}

\date{}

\title{Hamiltonicity of random graphs in the stochastic block model}
\author{Michael Anastos \\Freie Universit\"{a}t Berlin \\manastos@zedat.fu-berlin.de \and Alan Frieze\thanks{Research supported in part by NSF grant DMS1363136}\\Carnegie Mellon University\\alan@random.math.cmu.edu\and Pu Gao\thanks{Research supported by NSERC.}\\ University of Waterloo\\ pu.gao@uwaterloo.ca}
\begin{document}
\maketitle

\begin{abstract}
We study the Hamiltonicity of the following model of a random graph. Suppose that we partition $[n]$ into $V_1,V_2,\ldots,V_k$ and add edge $\{x,y\}$ to our graph with probability $p$ if there exists $i$ such that $x,y\in V_i$. Otherwise, we add the edge with probbability $q$. We denote this model by $\G(\bfn, p,q)$ and give tight results for Hamiltonicity, \jtb{including a critical window analysis,} under various conditions.
\end{abstract}

\section{Introduction}

The Hamiltonicity of various models of random graphs has been studied for many years. As far back as 1976, Koml\'os and Szemer\'edi announced their solution for the random graph $G_{n,m}$, although the published paper came out later in 1983 \cite{KS}. Since that time there have been many results on Hamiltonicity of random graphs, including but not restricted to,  binomial random graphs, random regular graphs,  binomial random graphs restricted to given minimum degrees, random $k$-out graphs, random percolation on given graphs, random graphs produced by (various types of) random graph processes, and also random hypergraphs.
See a recent bibliography~\cite{Fbib} by \jt{the second author} which goes into great details.

In this paper we study Hamiltonicity of random graphs from the so-called {\em Stochastic Block Model}. This random graph model has been the subject of much research in the computer science community. It is a generative model for social networks consisting of distinct communities. The model generalises the Erd\H{o}s-R\'{e}nyi random graphs, where every pair of vertices is connected by an edge independently with the same probability. In the stochastic block model, the probability of connecting a pair of vertices depends on which communities they belong to. Research on the stochastic block model is mainly on inferring the community membership given an instance sampled from the model. A recent paper by Abbe \cite{A} surveys this aspect.

A formal definition of the stochastic block model is given as follows.
Let $P$ be a symmetric $k\times k$ matrix with nonnegative entries between 0 and 1, and $\bfn=(n_1,\ldots,n_k)$ be a vector of positive integers. Let $n=\sum_{i=1}^n n_i$. Let $\G(\bfn,P)$ be a random graph constructed as follows. The vertex set is $V=\cup_{i=1}^k V_i$ where $V_i=\{(i, j), j\in [n_i]\}$, and any two vertices $(i_1,j_1)$ and $(i_2,j_2)$  are adjacent with probability $P(i_1,i_2)$ mutually independently. 
In this paper we consider the special case where $k\ge 1$ is a fixed integer, and $P$ has value $p=p(n)$ on the diagonal, and has value $q=q(n)$ off the diagonal. We denote $\G(\bfn, P)$ by $\G(\bfn, p,q)$ for this special $P$.

Unlike the Erd\H{o}s-R\'{e}nyi random graph, $\G(\bfn,p,q)$ is a non-homogeneous model where the distribution of the neighbourhood of vertex $v$ depends on which $V_i$ it belongs to. If $p=q$ then $\G(\bfn,p,q)$ reduces to $\G(n,p)$. If $p=0$ then $\G(\bfn,p,q)$ reduces to a random $k$-partite graph.
 The closest previous results to this work are the cases of Hamiltonicity of \jt{Erd\H{o}s-R\'{e}nyi graphs by Koml\'os and Szemer\'edi~\cite{KS}, and of} random bipartite graphs considered by \jt{the second author}\cite{Fbip}. The present paper utilises and extends the proofs in \jt{these papers} in a significant manner.

\section{The main results}
 We call vertex sets $V_i$ blocks, and an edge is called a {\em block edge} if its ends lie in the same block, and a {\em crossing edge} otherwise. Given a vertex $u\in V_i$, we say $u$ has {\em partition index} $i$. We aim to determine when $\G(\bfn, p,q)$ is Hamiltonian.

We will assume the following set of conditions.
\begin{align*}
\text{(A1)}& \quad \min_{1\le i\le k}\{pn_i+(n-n_i)q-\log n_i\}= \log\log n+O(1);\\
\text{(A2)}&  \quad qn^2=\omega(1);\\
\text{(A3)} & \quad \text{$\max_{1\le i\le k} n_i\le n/2$, if $p=O(1/n)$.}\\
\text{(A4)} & \quad \min_{1\le i\le k} n_i=\Omega(n).
\end{align*}
Note that if $\G(\bfn, p,q)$ is Hamiltonian then conditions (A2) and (A3) are necessary in general. If (A2) fails then with a non-vanishing probability there can be some $V_i$ such that $E(V_i, V\setminus V_i)=\emptyset$. If (A3) fails then $\G(\bfn,p,q)$ cannot be Hamiltonian if $p=0$. Condition (A4) can probably be relaxed, \jt{but it requires} more delicate analysis. \jtb{We will show that condition (A1) captures the critical window for Hamiltonicity of $\G(\bfn,p,q)$.}

Let $a_n$ and $b_n$ be two sequences of real numbers. We say $a_n=O(b_n)$ if there exists an absolute constant $C>0$ such that $|a_n|\le C|b_n|$ for every $n\ge 1$. We say $a_n=o(b_n)$ if $b_n>0$ for all $n\ge 1$ and $\lim_{n\to\infty} a_n/b_n=0$. If $a_n>0$ for all $n\ge 1$, and $b_n=O(a_n)$ (or $b_n=o(a_n)$) then we write $a_n=\Omega(b_n)$ (or $b_n=\omega(a_n)$ respectively). \jt{We will consider a sequence of random graphs indexed by their order, denoted by $n$, and all constants in this paper do not depend on $n$. All asymptotics refer to $n\to\infty$.} Given a graph property $\Gamma$, we say $\G(n,p,q)\in \Gamma$ asymptotically almost surely (a.a.s.) if $\lim_{n\to \infty}(\G(n,p,q)\in \Gamma)=1$. 
Let \ham\ denote the class of Hamiltonian graphs. Our main result is the following.

\begin{thm}\lab{thm:ham}
Assume $p$ and $q$ and $\bfn$ satisfy assumptions (A1)--(A4). 
\[
\lim_{n\to\infty}\pr(\G(\bfn,p,q)\in \ham)=\exp\left(-\sum_{1\le i\le k} e^{-c_i}\right),
\]
where $c_i=pn_i+(n-n_i)q-\log n_i-\log\log n$.
\end{thm}

As Hamiltonicity is an increasing property, the following corollary follows immediately.
\begin{cor}\lab{cor:ham}
Assume $p$ and $q$ and $\bfn$ satisfy assumptions (A2)--(A4). 
\[
\lim_{n\to\infty}\pr(\G(\bfn,p,q)\in \ham)=\left\{
\begin{array}{ll}
0& \mbox{if $\min_{1\le i\le k}\{pn_i+(n-n_i)q- \log n_i\}<\log\log n  -\omega(1)$}\\
1 & \mbox{if $\min_{1\le i\le k}\{pn_i+(n-n_i)q- \log n_i\}> \log\log n +\omega(1)$}.
\end{array}
\right.
\]
\end{cor}

\section{Small degrees}

Let \degree\ denote the class of graphs with minimum degree at least 2. Note that $G\notin \degree$ implies that $G\notin\ham$. Thus, the following lemma immediately yields an upper bound on the probability that $\G(\bfn,p,q)$ is Hamiltonian.

\begin{lemma}\lab{lem:mindegree}
Assume (A1) \jt{and (A4)}. Then,
\begin{enumerate}
\item[(a)] \[
\lim_{n\to\infty}\pr(\G(\bfn,p,q)\in\degree)=\exp\left(-\sum_{1\le i\le k} e^{-c_i}\right),
\]
where $c_i=pn_i+(n-n_i)q-\log n_i-\log\log n$.
\item[(b)] For constants $0<\alpha<1$, and $\beta>\alpha+\alpha\ln(1/\alpha)$. A.a.s.\ $\G(\bfn,p,q)$ contains at most $n^{\beta}$ vertices whose degree is at most $\alpha \log n$.
\end{enumerate}
\end{lemma}
\proof For part (a), let $X_j(i)=\sum_{v\in V_i} \ind{d(v)=j}$ be the number of vertices in $V_i$ with degree $j$. Let $W_1$ and $W_2$ be two independent random variables with $W_1\sim \Bin(n_i-1,p)$ and $W_2\sim\Bin(n-n_i,q)$. Let $j=O(\log n)$. By (A1) and (A4) we have $p^2n=o(1)$, $pj=o(1)$, $j^2=o(n_i)$, and $j^2=o(n-n_i)$. Then,
\begin{align}
\ex X_j(i) &= n_i \pr(W_1+W_2=j)\non\\
&=n_i \sum_{s=0}^j \binom{n_i-1}{s} p^{s}(1-p)^{n_i-1-s} \binom{n-n_i}{j-s}q^{j-s}(1-q)^{n-n_i-j+s}\non\\
&\sim n_i \exp(-pn_i-q(n-n_i))\sum_{s=0}^j \frac{n_i^s(n-n_i)^{j-s}}{s!(j-s)!}p^sq^{j-s}\non\\
&=n_i  e^{-\phi_i}\frac{\phi_i^j}{j!}, \lab{degree-expectation}
\end{align}
where $\phi_i=pn_i+q(n-n_i)$. By (A1), $\phi_i-(\log n_i+\log\log n)>C$ for some constant real $C$ and for all $i\in[k]$.
It follows immediately that a.a.s.\ $X_0(i)=0$ for every $i$.

Recall that $c_i=\phi_i-(\log n_i+\log\log n)$.  Let $X_1=\sum_{i\in [k]} X_1(i)$. Then, $\ex X_1\sim \sum_{i\in [k]} e^{-c_i}$. By (A1), $\ex X_1=\Theta(1)$.
Using the standard method of moments (we omit the tedious calculations), it is easy to prove that $X_1$ is asymptotically Poisson. Hence, 
\be\lab{degree1}
\pr(X_1=0)\sim \exp\left(-\sum_{i\in [k]}e^{-c_i}\right).
\ee
The lemma follows by~\eqn{degree1} and the fact that a.a.s.\ $X_0(i)=0$ for every $i\in[k]$.

For part (b), from~\eqn{degree-expectation} we have
\bean
\sum_{j\le \alpha \log n}\sum_{i\in [k]} \ex X_j(i) &<& (1+o(1)) \sum_{j\le \alpha \log n}\sum_{i\in [k]}  n_i \exp(-\log n_i) \frac{(\log n_i)^j}{j!}\\
&<& k \sum_{j\le \alpha \log n} \left(\frac{e\log n}{j}\right)^j=(1+o(1))k n^{\rho(\alpha)},
\eean
where $\rho(\alpha)=\alpha+\alpha\log(1/\alpha)$. Part (b) follows by the Markov inequality.
\qed
\section{Vertex expansion and connectivity}

Let $G$ be a graph and $S\subseteq V(G)$, define
\bean
\calN_G(S)&=&\{j\in V(G)\setminus S:\ \exists i\in S, i\sim j\}\\
N_G(S)&=&|\calN_G(S)|\\
n_1(G)&=&\sum_{i\in V(G)} \ind{d(i)\le 1}.
\eean
I.e.\ $\calN_G(S)$ is the set of vertices not in $S$ which are adjacent to some vertex in $S$ in graph $G$.  We may drop $G$ from the subscript if the underlying graph $G$ is clear from the context. 
\begin{definition}
A vertex in $G$ is called  small, if its degree is less than $\log n/10$. A vertex with degree at least $\log n/10$ is called large.
\end{definition}
We say $G$ has property \expand, if: there exists a constant $\eps_0>0$ such that
\[
\mbox{for every $S\subseteq V(G)$ where $|S|\le \eps_0 n$, $|N_G(S)|\ge 2|S|\cdot \ind{n_1(G)=0}$.}
\]
If $F$ is a subset of edges in $G$, we use $G-F$ to denote the subgraph of $G$ obtained by deleting edges in $F$.
 We say $G$ has property \expandd, if  the following holds:
For any $F\subseteq E(G)$, such that $|F\cap \calN_G(v)|=0$ if $v$ is small, and is at most $\log n/100$ if $v$ is large, \jt{we have that}  $G-F$ is connected and $G-F\in \expand$.

\begin{lemma}\lab{lem:expansion} 
Assume (A1), (A2) and (A4). Then a.a.s. $\G(\bfn,p,q)\in\expandd$.
\end{lemma}
Its  technical proof is postponed till Section~\ref{sec:technical}.

\section{\jt{Overview of the proof} of Theorem~\ref{thm:ham}}

Since $G\notin \degree$ implies that $G\notin \ham$, the upper bound of $\pr(\G(\bfn,p,q)\in\ham)$ is implied by Lemma~\ref{lem:mindegree}(a). Next we prove the lower bound. Let \length\ denote the set of graphs where a longest path contains the same number of vertices as in a longest cycle, and let \connect\ denote the class of connected graphs. Note that if $G\in\length\cap\connect$ then $G$ must be Hamiltonian, since otherwise, by connectivity it is always possible to extend a longest cycle into a path which contains more vertices than the cycle we start with, contradicting \jt{with} $G\in\length$. It follows then that
\be
\pr(\G(\bfn,p,q)\in\ham)\ge \pr(\G(\bfn,p,q)\in\length\cap\connect)\ge \pr(\G(\bfn,p,q)\in\length\cap\connect\cap\degree).\lab{lowerbound}
\ee
Our goal is to prove that $\pr(\G(\bfn,p,q)\in\length\cap\connect\cap\degree)\sim \pr(\G(\bfn,p,q)\in\degree)$, which then yields the asymptotic probability desired by Theorem~\ref{thm:ham}. The proof of the lower bound of $\pr(\G(\bfn,p,q)\in\ham)$ will be split into three cases: (1) $p,q=\omega(1/n)$; (2) $p=O(1/n)$; and (3) $q=O(1/n)$. In all three  cases, we will use a multi-round exposure technique of $\G(\bfn,p,q)$. Roughly speaking, we will expose a subgraph $G_b\subseteq G$ where $G\sim \G(\bfn,p,q)$ and $G_b$ contains most edges of $G$.   Case 1 is the simplest case, in which  we will define graph property \typical\ which consists of a set of properties that hold a.a.s.\ for $\G(\bfn,p,q)$.
Then  we will define \col\ to be a set of properties that edges in $G\setminus G_b$ must satisfy. Then, we will prove that 
\[
\pr(\{G\in \overline{\length}\cap \connect\cap\degree\cap\typical\}\cap \col)\ll \pr(\col \mid G\in  \overline{\length}\cap \connect\cap\degree\cap\typical).
\]
This implies that $\pr(G\in \overline{\length}\cap \connect\cap\degree\cap\typical)=o(1)$, which will lead us to derive the asymptotic probability for $\pr(G\in\length\cap\connect\cap\degree)$. While \jt{obtaining a lower bound for} $\pr(\col \mid G\in  \overline{\length}\cap \connect\cap\degree\cap\typical)$ is rather straightforward, an upper bound for $\pr(\{G\in \overline{\length}\cap \connect\cap\degree\cap\typical\}\cap \col)$ is obtained by using  {\em P\'osa rotations} and bounding the probability that the longest path does not get extended by the edges exposed in the second stage. This is a standard technique for proving Hamiltonicity in random graphs.

In Case 2, $\G(\bfn,p,q)$ is similar to the random $k$-partite graph $\G(\bfn,0,q)$. Two complications arise in this case. Firstly, we cannot totally ignore block edges in $\G(\bfn,p,q)$ as they contribute to degree 2 vertices in $\G(\bfn,p,q)$ with a non-varnishing probability. More specifically, $\pr(\G(\bfn,p,q)\in \ham)>\pr(\G(\bfn,0,q)\in\ham)$ when $p=\Theta(1/n)$ and thus, the proof cannot be obtained by simply extending the proof for random bipartite graphs to random $k$-partite graphs. Secondly, due to the asymmetry between $p$ and $q$, the edges exposed in later stages will not be uniformly distributed and we need to take care of the multipartition of the vertices. This is similar to the case of the random bipartite graphs. 

In Case 3, $\G(\bfn,p,q)$ looks like a collection of $G_i\sim \G(n,p)$ plus a set of random edges between every pair of $G_i$, $G_j$, $1\le i<j\le k$. A tempting approach would be to find a Hamilton cycle in each $G_i$ and then somehow connect these cycles by using a few crossing edges to form a Hamilton cycle in $\G(\bfn,p,q)$. This approach would succeed if $q=o(1/n)$. However, when $q=\Theta(1/n)$, similar to Case 2, the crossing edges are contributing, with a non-negligible probability, to the degree 2 vertices in $\G(\bfn,p,q)$. Thus, we cannot purely focus on structures in $G_i$. Instead, inside each $G_i$, \jtb{we will take particular care of the vertices with degree less than 2, and} we will look for a small number of vertex disjoint paths covering all vertices in $G_i$. These paths have specified end vertices. Then we will stitch these paths together with some crossing edges to form a Hamilton cycle in $\G(\bfn,p,q)$.

As part of the overview of the proof, we define \typical\ and \col. They will be used in the proof of the first case, and in the second case as well with some minor modifications. In Case 3, their definitions will be significantly modified.

\subsection{\typical}

We say $G\in\typical$ if $G$ satisfies the following set of properties.
\begin{enumerate}
\item[(T1)] $G\in \expandd$.
\item[(T2)] There are at most $n^{0.4}$ small vertices.
\item[(T3)] If $p=o(\log n/n)$ then every vertex is incident with at most $\log n/200$ block edges. If $q=o(\log n/n)$ then every vertex is incident with at most $\log n/200$ crossing edges.
\item[(T4)] The maximum degree is at most $C\log n$ for some sufficiently large constant $C>0$.
\end{enumerate}

\begin{lemma} \lab{lem:typical}
A.a.s.\ $G\in\typical$.
\end{lemma}
\proof (T1) follows by Lemma~\ref{lem:expansion}. \jtb{(T2) follows by Lemma~\ref{lem:mindegree}(b).} \jt{(T3) and} (T4) follow by a standard first moment argument \jtb{similar to the proof of Lemma~\ref{lem:mindegree}}. We omit the details.\qed

\subsection{\col}
\lab{sec:col}
Let $L(G)$ denote the length of a longest path in $G$.
Assume $G'\subseteq G$. Let $F=E(G\setminus G')$.
We say $(G,G')\in\col$ if 
\begin{enumerate}
\item[(a)] $|F\cap \calN_{G}(v)|$ is 0 if $v$ is small in $G$, and is at most $\log n/100$ if $v$ is large in $G$;
\item[(b)] $L(G)=L(G')$ if $L(G')<n-1$, and $G\notin \ham$ if $L(G')=n-1$.
\end{enumerate}
In the proof of Theorem~\ref{thm:ham}, we will use \col\ to denote the event that $(G,G_b)\in\col$, although $G$ and $G_b$ are defined differently in the three cases. We will recall the definition of \col\ when we proceed to the proof in each case.

\section{Proof of Theorem~\ref{thm:ham}: when $p,q=\omega(1/n)$}
\lab{sec:case1}

We first define $G_b$.
\subsection{$G_b$}
Let $\barp=a/n\log n$, where \jt{$a=1$. We will choose similar parameters for multi-round exposures in Cases 2 and 3 with different values of $a$. We keep $a$ in the definition of $\barp$ for the ease of comparison.} Define 
\[
p_1=1-\frac{1-p}{1-\barp};\quad q_1=1-\frac{1-q}{1-\barp}.
\]
In case 1, both $p_1$ and $q_1$ are real numbers between 0 and 1. We will run a two stage exposure  of the edges in $\G(\bfn,p,q)$. First, generate $G_b\sim \G(\bfn,p_1,q_1)$, then independently for every non-edge $x$ in $G_b$, add $x$ to the graph with probability $\barp$. Call the resulting graph $G$. It is straightforward to verify that $G\sim\G(\bfn,p,q)$. For convenience, colour the edges in $G_b$ blue and the edges in $E(G\setminus G_b)$ red.

By our definition of $p_1$ and $q_1$ it is easy to see that (A1) and (A2) are satisfied with $p$ and $q$ replaced by $p_1$ and $q_1$. Recall that $\col$ denotes the event that $(G,G_b)\in\col$.

The next two lemmas bound $\pr(\col\mid G\in \overline{\length}\cap \degree\cap\typical)$ and $\pr(\{G\in \overline{\length}\cap\degree\cap\typical\}\cap \col)$.
\begin{lemma}\lab{lem:colour}
There exists a function $f=f(n)\to\infty$ as $n\to \infty$ such that
\[
\pr(\col\mid G\in \overline{\length}\cap \degree\cap\typical)\ge \exp(-an/f\log n).
\]
\end{lemma}

\begin{lemma}\lab{lem:length}
\[
\pr(\{G\in \overline{\length}\cap\degree\cap\typical\}\cap \col) \le \exp(-\Omega(an/\log n)).
\]
\end{lemma}

Now we are ready to prove Theorem~\ref{thm:ham} in case 1.

\no{\em Proof of Theorem~\ref{thm:ham} (case 1).} 
By Lemmas~\ref{lem:colour} and~\ref{lem:length}, 
\[
\pr(G\in \overline{\length}\cap \degree\cap\typical)\le \frac{\exp(-\Omega(an/\log n))}{\exp(-an/f\log n)}.
\]
As $f\to\infty$ as $n\to\infty$, the above probability is $o(1)$. By Lemmas~\ref{lem:expansion} and~\ref{lem:typical}, $\pr(G\in \connect\cap\typical )=1-o(1)$. It follows immediately that 
\[
\pr(G\in \length\cap\degree\cap\connect)=\pr(G\in\degree)-\pr(G\in \overline{\length}\cap \connect\cap\degree\cap\typical)+o(1)=\pr(\degree)+o(1).
\]
 By~\eqn{lowerbound} and the fact that $G\in\ham$ implies $G\in \degree$, we have $\pr(G\in\ham)= \pr(G\in\degree)+o(1)$.  Together with Lemma~\ref{lem:mindegree} this yields the asymptotic probability of $\pr(G\in\ham)$ as in Theorem~\ref{thm:ham}.\qed

It remains to prove Lemmas~\ref{lem:colour} and~\ref{lem:length}.

\subsection{Proof of Lemma~\ref{lem:colour}}

Equivalently we can define $G_b$ as follows. Take $G\sim\G(\bfn,p,q)$. Define
\[
p^*=\frac{\bar p(1-p)}{(1-\bar p)p},\quad q^*=\frac{\bar p(1-q)}{(1-\bar p)q}.
\]
Do the following independently for every edge $x\in G$: if $x$ is a block edge, delete $x$ with probability $p^*$; if $x$ is a crossing edge, delete $x$ with probability $q^*$. As $p(1-p^*)=p_1$ and $q(1-q^*)=q_1$ with our definition of $(p^*,q^*)$, we immediately have
\begin{claim}
The resulting graph is distributed as $\G(\bfn,p_1,q_1)$.
\end{claim}

We will prove that conditioning on $G=H$ for any $H\in \overline{\length}\cap \degree\cap\typical$, 
$\pr(\col\mid G=H)\ge \exp(-an/f\log n)$, and Lemma~\ref{lem:colour} follows. 

Consider the set of edges deleted in generating $G_b$ from $H$. Colour these edges red.

Let $P$ be a longest path in $H$. Note that \col\ is implied if 
\begin{enumerate}
\item[(B1)] no large vertex in $H$ is incident with more than $\log n/100$ red edges.
\item[(B2)] no small vertex in $H$ is incident with a red edge;
\item[(B3)] no edge in $P$ is red. 
\end{enumerate}

Let $\calX$ be the union of the set of edges in $P$ and the set of edges incident with small vertices. Let $X=|\calX|$. By (T2), $|\calX|\le (n-1)+n^{0.4} \log n/10<2n$. As $\max\{p^*,q^*\}=o(a/\log n)$ in case 1, the probability that none of the edges in $\calX$ is deleted is
at least $(1-o(a/\log n))^{2n}\ge \exp(-an/f\log n)$ for some $f\to\infty$. Note that (B2) and (B3) are implied if no edges in $\calX$ are red. Hence, $\pr(B1\cap B_2)\ge \exp(-an/f\log n)$. 

Let $\overline{\calX}$ be the set of edges in $H$ that are not in $\calX$. Condition on no edges in $\calX$ were deleted (i.e.\ became red). 
We will prove that a.a.s.\ every vertex is incident with at most $\log n/200$ red block edges in $\overline{\calX}$, and at most $\log n/200$ red crossing edges in $\overline{\calX}$. By (T3) we may assume that $p,q=\Omega(\log n/n)$. By (T4), each vertex has degree $\jt{O(\log n)}$. By the definition of $p^*$ and $q^*$, every edge is deleted (i.e.\ becomes red) with probability $O(\barp\cdot\max\{1/p,1/q\})=O(1/\log^2 n)$.
By the tail bounds for the binomial distribution and the union bound, a.a.s.\ every vertex is incident with at most $o(\log n)$ red edges and thus $\pr(B1\mid B2\cap B3)=1-o(1)$. 
 Hence, $\pr(\col\mid G=H)\ge \pr(B1\cap B2\cap B3)\ge \exp(-an/f\log n)$ for some $f\to\infty$. Lemma~\ref{lem:colour} follows. \qed

\subsection{Proof of Lemma~\ref{lem:length}}

Recall the definition of \expand\ and \expandd. 
Assuming (T1), there exists an absolute constant $\eps_0>0$ such that
\begin{equation}
\text{for every $S\subseteq V(G)$ where $|S|\le \eps_0 n$, we have $|N_{G}(S)|\ge 2|S|\cdot\ind{n_1(G)=0}$.} \lab{eps1}
\end{equation}
 By the definition of \col\ we immediately have the following claim.
\begin{claim}\lab{claim:imply}
$\{G\in \overline{\length} \cap\degree\cap\typical\}\cap\col$ implies $B\cap\{G_b\in \connect\cap\expand\cap\degree\}$, where
$B=\{L(G)=L(G_b)<n-1\}\cup (\{L(G_b)=n-1\} \cap\{ G\notin\ham\})$.
\end{claim}
\jt{
\proof If $G\in\typical\cap\degree$ and $(G,G_b)\in\col$, then by (T1) we have $G_b\in  \connect\cap\expand\cap\degree$. Moreover, $(G,G_b)\in\col$ implies $B$. This proves our claim.\qed\smallskip
}

Hence, it is sufficient to prove
\[
\pr(B\mid G_b\in \connect\cap\expand\cap\degree)\le \exp(-\Omega(an/\log n)),
\]
as
\bean
\pr(\{G\in \overline{\length} \cap\degree\cap\typical\}\cap\col)&\le& \pr (B\cap\{G_b\in \connect\cap\expand\cap\degree\})\\
&\le& \pr(B\mid G_b\in \connect\cap\expand\cap\degree),
\eean
by Claim~\ref{claim:imply}.
Note that $G$ is obtained by adding every non-edge in $G_b$ independently with probability $\barp$. We will prove that conditioning on any graph $G_b\in\connect\cap\expand$, adding approximately $\bar p\binom{n}{2}\sim an/2\log n$ edges will either increase $L(G_b)$, or complete a Hamilton path in $G_b$ to a Hamilton cycle, with sufficiently high probability. We will use the classical technique of P\'osa rotations to bound the probability of $B$.

{\em P\'osa rotations.} Let $P=v_0,v_1,\ldots, v_{\ell}$ be a longest path in $G_b$. Then $v_0$ is not adjacent to $v_{\ell}$, and all the neighbours of $v_{\ell}$ in $G_b$ must be in $P$, since otherwise we can extend $P$ to a longer path. Assume $v_iv_{\ell}$ is an edge in $G_b$ where $i<\ell-1$, then another longest path $P'=v_0,\ldots,v_iv_{\ell},v_{\ell-1},\ldots,v_{i+1}$ can be obtained by using the edge $v_iv_{\ell}$ instead of $v_iv_{i+1}$. This operation from $P$ to $P'$ is called a P\'osa rotation. Consider the set $\sP$ of longest paths obtained by repeatedly rotating $P$. All of these paths start from $v_0$ and end at a vertex that is in $P$. Let $\text{End}(v_0)$ denote the set of ends  other than $v_0$ in the paths in $\sP$.
A key observation is the following. The reader may refer to~\cite{FK} for a proof. 
\begin{lemma}
$|N_{G_b}(\text{End}(v_0))|<2|\text{End}(v_0)|$.
\end{lemma}
As $G_b\in\expand\cap\degree$ we immediately have that $|\text{End}(v_0)|\ge \eps_0 n$ where $\eps_0$ is specified in~\eqn{eps1}. 

Now for every $v\in \text{End}(v_0)$, there is a longest path $P_v$ which is obtained from $P$ by repeatedly applying P\'osa rotations. Let $\text{End}(v)$ denote the set of ends other than $v$ in the longest paths obtained by rotating $P_v$. Again, we have
$|N_{G_b}(\text{End}(v))|<2|\text{End}(v)|$, which implies that $|\text{End}(v)|\ge \eps_0 n$. Consider the set $\sE$ of pairs of vertices $(x,y)$ where $x\in \text{End}(v_0)$, and $y\in \text{End}(x)$. We have that $|\sE|\ge \eps_0^2 n^2$. Moreover, adding any pair in $\sE$ as an edge to $G$ will either form a Hamilton cycle in $G$, if $\ell=n-1$, or form a cycle with length $\ell+1$, and then using the fact that $G$ is connected, we can extend the cycle to a path of length $\ell+1$, if $\ell<n-1$. \jt{In either case, event $B$ fails.} For that reason, we call $\sE$ a set of {\em boosters}. We have shown that $B$ fails if $\sE\cap E(G\setminus G_b)\neq \emptyset$.

By the construction of $G$, every edge in $\sE$ is added to $G$ in the second stage of edge exposure, independently with probability $\barp$. The probability that none of these edges are added is $(1-\barp)^{|\sE|/2}\le \exp(-\eps_0^2 a n/2\log n)$. This completes the proof for Lemma~\ref{lem:length}. \qed

\section{Proof of Theorem~\ref{thm:ham}: when $p=O(1/n)$}

In this case, we will define a 3-round edge exposure of $\G(\bfn,p,q)$.
\subsection{$G_b$ and $G_y$} 
Let $\bar q=a/n\log n$ where $a\to \infty$ and $a=o(\log n)$. Define
\[
q_1=1-\frac{1-q}{(1-\bar q)^2}.
\] 

First, generate $\hat{G}\sim \G(\bfn,0,q_1)$. 
Then, expose block edges that are incident to vertices with degree at most one in $\hat{G}$.  Let $G_b$ be the resulting graph. All edges in $G_b$ are coloured blue. By bounding the expected number of vertices having degree at most 1 in $\hat{G}$  by $o(\log n)$, we immediately have the following claim.
\begin{claim}\lab{claim:blockedge}
A.a.s. fewer than $\log n$ vertices in $G_b$ are incident to a block edge.
\end{claim}

 Then, for every crossing non-edge $x$ in $G_b$, add $x$ to the graph $G_b$ with probability $\bar q$. All edges added in the second stage are coloured yellow. Call the resulting graph $G_y$. Finally, for each crossing non-edge in $G_y$, add it to the graph with probability $\bar q$ and colour it red.  The crossing edges of the resulting graph $G$ have the same distribution as those in $\G(\bfn,p,q)$. Obviously we can couple $G$ with $G'\jt{\sim}\G(\bfn,p,q)$ such that $G\subseteq G'$. It is thus sufficient to prove that the lower bound for $\pr(G\in\ham)$ matches the lower bound in Theorem~\ref{thm:ham}.

Let now \typical\ include properties (T1)--(T4) as well as the following property.
\[
(T5): \quad \mbox{Fewer than $\log n$ vertices of $G$ are incident with fewer than 2 crossing edges.}
\]

Let \col\ denote $(G,G_b)\in \col$ as defined in Section~\ref{sec:col}. Note that (T3), (T5) and \col\ implies the following.
\[
\col1: \quad \mbox{Fewer than $\log n$ vertices in $G_b$ are incident to a block edge.}
\]

The next two lemmas bound $\pr(\col\mid G\in \overline{\length}\cap \degree\cap\typical)$ and $\pr(\{G\in \overline{\length}\cap\degree\cap\typical\}\cap \col)$ in case 2.
\begin{lemma}\lab{lem:colour2}
There exists a constant $K>0$ such that
\[
\pr(\col\mid G\in \overline{\length}\cap \degree\cap\typical)\ge \exp(-Kan/\log^2 n).
\]
\end{lemma}

\begin{lemma}\lab{lem:length2}
\[
\pr(\{G\in \overline{\length}\cap\degree\cap\typical\}\cap \col) \le \exp(-\Omega(a^2n/\log^2 n)).
\]
\end{lemma}

Since $\pr(\{G\in \overline{\length}\cap\degree\cap\typical\}\cap \col) \ll \pr(\col\mid G\in \overline{\length}\cap \degree\cap\typical)$
as $a\to \infty$, the proof of Theorem~\ref{thm:ham} for case 2 follows exactly as in case 1. It only remains to prove Lemmas~\ref{lem:colour2} and~\ref{lem:length2}. 

\subsection{Proof of Lemma~\ref{lem:colour2}}
The proof is basically the same as that for Lemma~\ref{lem:colour}. The probability bound is only different because of the different range for $p$ and $q$. We only point out the differences in the proof. We give an equivalent definition of $(G_b,G)$ as follows. Take $G'\sim \G(\bfn,p,q)$. Let $q^*$ be such that $q(1-q^*)=q_1$. Independently delete each crossing edge $x$ in $G'$ with probability $q^*$. Call the resulting graph $G_b'$. Next delete all block edges in $G_b'$ except for the edges incident with a vertex which is only incident with at most one crossing edge in $G_b'$. The resulting graph has the same distribution as $G_b$. Let $G$ be the graph obtained by taking all crossing edges in $G'$ and all block edges in $G_b$. Note that $G$ and $G_b$ only differ on the crossing edges.

Now $\max\{p^*,q^*\}$ in the proof of Lemma~\ref{lem:colour} is replaced by $q^*$. By (A1), (A5) and the assumption $p=O(1/n)$ we have $q=\Theta(\log n/n)$ and thus $q^*=\Theta(a/\log^2 n)$. It follows that the probability that none of the edges in  $\chi$ is deleted is at least $(1-\Theta(a/\log^2 n))^{2n}\ge \exp(-\Theta(an/\log^2 n))$. The rest of the proof is the same as in Lemma~\ref{lem:colour}. \qed

\subsection{Proof of Lemma~\ref{lem:length2}}
Recall the definition of \typical, \expand\ and \expandd. Assuming $G\in\expandd$, there exists an absolute constant $\eps_0>0$ such that~\eqn{eps1} holds. Claim~\ref{claim:imply} continues to hold. 
Thus, it is sufficient to show
\[
\pr(B\mid\{G_b\in \connect\cap\expand\cap\degree\})\le \exp(-\Omega(a^2n/\log^2 n)).
\]
As in the proof of Lemma~\ref{lem:length}, we will show that the probability that the additional yellow and red edges do not increase $L(G_b)$ is small. 
However, the graph induced by the yellow and red edges is $k$-partite.
and so
there is a subtle issue here with $k$-partite graphs. The yellow and red edges exposed in the second and third stage must respect to the vertex partition. If the set of boosters, the potential non-edges whose addition will allow an extension of the longest path, are all unichromatic, i.e. the ends of each booster have the same partition index, then the edges exposed in the second stage will not help to extend the paths. The purpose of using 3 rounds of edge exposure is to cope with this vertex partition issue. To explain how it works, we need a few definitions.

Let $P=(v_0,\ldots,v_{\ell})$ be a longest path in the graph $G_b$. Let
$\sP$ be a list of longest paths obtained as follows. $\sP=\{P\}$ initially and $\text{End}(v_0)=\{v_{\ell}\}$. Take a path $P'\in \sP$ and do a P\'osa rotation, if it creates a path whose end other than $v_0$ is not yet in $\text{End}(v_0)$. Then add this new path to $\sP$ and add its end vertex other than $v_0$ to $\text{End}(v_0)$. Since $|\text{End}(v_0)|<n$, this process will terminate. As before, given that $G_b\in\expand$ we have  $|\text{End}(v_0)|\ge \eps_0 n$.

For each $v$, let $\sigma(v)$ denotes its partition index.
For each $x\in \text{End}(v_0)$, let $P_x$ be the path that was added to $\sP$ when $x$ was added to $\text{End}(v_0)$. Let $\phi(x)$ be the vertex adjacent to $x$ on $P_x$. Let $\sP(x)$ denote the set of longest paths obtained by rotating $P_x$ with $x$ being the fixed end, and let $\text{End}(x)$ denote the set of the ends of the paths other than $x$ in $\sP(x)$. Then, all paths in $\sP(x)$ must start with $x\phi(x)$. Moreover, $\sigma(x)\neq\sigma(\phi(x))$, unless $x$ is incident with a block edge.

We call a vertex $x$ in $\text{End}(v_0)$ {\em good}, if at least $\frac{1}{2}|\text{End}(x)|$ vertices in $\text{End}(x)$ have distinct partition index from $x$.

Consider $G_y$. Let $\tt{Out}$ denote the set of vertices $x$ for which there exists a longest path $P$ in $G_y$ such that $x\notin P$ and one end of $P$ lies in distinct vertex part from $x$. 
We will show the following.
\begin{claim}\label{c:out}
Either there are $\frac{1}{2}|\text{End}(v_0)|$ good vertices in $\text{End}(v_0)$, or with probability $1-\exp(-\Omega(an/\log n))$, either $L(G_y)>L(G_b)$ or $|{\tt Out}|=\Omega(an/\log n)$ in $G_y$.
\end{claim}
Note $\exp(-\Omega(an/\log n))\le \exp(-\Omega(a^2n/\log^2 n))$ as $a\to \infty$ is chosen to be $o(\log n)$.

If there are $\frac{1}{2}|\text{End}(v_0)|$ good vertices in $\text{End}(v_0)$, then with the same proof as in case 1, we can bound $\pr(B\mid G_b\in\connect\cap \expand\cap \degree)$ by $\exp(-\Omega(an/\log n))\le \exp(-\Omega(a^2 n/\log^2 n))$. 

If instead $|{\tt Out}|=\Omega(an/\log n)$ in $G_y$, then we will prove the following claim which, together with Claim~\ref{c:out} and the argument above, completes the proof for Lemma~\ref{lem:length2}.

\begin{claim}
If $|{\tt Out}|=\Omega(an/\log n)$ then $B$ holds with probability $\exp(-\Omega(a^2n/\log^2 n))$.
\end{claim}
\proof Let $x\in \tt{Out}$ and let $P$ be the corresponding path, and $y$ be the end on $P$ where $\sigma(x)\neq \sigma(y)$. Consider P\'osa rotations on $P$ with $y$ as the fixed end, and let $\text{End}(y)$ denote the ends obtained. Either at least half of $\text{End}(y)$ have partition indices distinct from $\sigma(x)$, or at least half with partition indices distinct from $\sigma(y)$. This implies that $|\sE_x|=\Omega(n)$ where 
\[
\sE_x=\{\{x,y\}:\ z\in \text{End}(y), \sigma(x)\neq \sigma(z)\}\cup \{\{y,z\}:\ z\in \text{End}(y), \sigma(y)\neq \sigma(z) \}.
\] 
Let $\sE=\cup_{x\in {\tt Out}}\sE_x$. Then $|\sE|=\Omega(a n^2/\log n)$. Moreover, $B$ fails if $\sE\cap E(G\setminus G_y)\neq\emptyset$. Since every booster in $\sE$ appears in the second stage of edge exposure with probability $\bar q$. The probability that none of them appears is at most $(1-\bar q)^{|\sE|}=\exp(-\Omega(a^2 n/\log^2 n))$. \qed\ss

\no {\em Proof of Claim~\ref{c:out}.\ } Let $\text{Bad}(v_0)$ denote the set of vertices in $\text{End}(v_0)$ that are not good and assume $|\text{Bad}(v_0)|>\frac{1}{2}|\text{End}(v_0)|$. Let $x\in \text{Bad}(v_0)$. Let $\tt{Mono}(x)$ denote the set of vertices in $\text{End}(x)$ which have the same partition index as $x$. By the definition of $\text{Bad}(v_0)$ and because $x\in \text{Bad}(v_0)$, $|\tt{Mono}(x)|\ge \frac12 |\text{End}(x)|$. Consequently, the vertices in $\tt{Mono}(x)$ must have distinct partition indices from $\phi(x)$, unless $x$ is incident with a block edge in $G_b$. Let $\tt{A}$ denote the set of vertices that are not incident with any block edges in $G_b$.  Consider all pairs of vertices $\sE'=\{(\phi(x),y): x\in \text{Bad}(v_0)\cap{\tt A},y\in \tt{Mono}(x)\}$. As $\typical\cap\col$ implies \col1, it follows that $|\sE'|=\Omega(n^2)$. If any yellow edge exposed in the second round of edge exposure is in $\sE'$, then we find a cycle $C$ by deleting $x$ from the path $P'$, and then adding the edge $y\phi(x)$. $G_b$ is connected and by (A3) there must exist a vertex $z$ such that (a) $z$ has distinct partition index from $x$; (b) $z$ is not on $P'$; (c) there is a path from $z$ to $P'$. Using the cycle $C$ and the path from $z$ to $P'$ we obtain a path of length greater than that \jt{of} $P'$, which implies $L(G_y)>L(G_b)$, or we obtain a path of the same length as $P'$ and one of its ends, namely $z$, has a distinct partition index from $x$. This implies that $x\in \tt{Out}$. Assume it is always the latter case for each yellow edge in $\sE'$. Since $|\sE'|=\Theta(n^2)$, it follows that $\ex|{\tt Out}|=\bar q \cdot \Theta(n^2)=\Omega(an/\log n)$. The probability bound in the claim follows by the Chernoff bound.\qed

\section{Proof of Theorem~\ref{thm:ham}: when $q=O(1/n)$.} 

In this case we have $\min_{1\le i\le k}\{pn_i- \log n_i\}=\log\log n+O(1)$. The subgraph induced by the vertices in block $i$ 
is distributed as $\G(n_i,p)$. 

Let $\bar p=a/n\log n$ \jt{with $a=1$}. Define
\[
p_1=1-\frac{1-p}{(1-\bar p)}.
\] 

We give a quick overview of the proof in this case. First we generate $G_b\sim \G(\bfn, p_1,q)$. We call a vertex ``problematic'' if it is incident with fewer than 2 block edges in $G_b$. Let ${\mathcal P}$ denote the set of problematic vertices.  All edges in $G_b$ are coloured blue. In the second round of edge exposure, we expose block edges that are not present in $G_b$ with probability $\bar p$. The resulting graph is $G$.  All edges in $G\setminus G_b$ are coloured red. It is easy to see that $G\sim \G(\bfn,p,q)$. We obtain our Hamilton cycle by (i) finding a collection ${\mathcal A}_1$ of paths of length 1 or 2 that cover the problematic vertices, \jtb{(ii) finding a collection ${\mathcal A}_2$ of crossing edges so that the number of crossing edges in ${\mathcal A}_1\cup {\mathcal A}_2$ between any pairs $V_i$ and $V_j$ is even and positive, and (iii)} finding a collection of vertex disjoint paths that connect the ends of the paths \jtb{and edges} in ${\mathcal A}\jtb{={\mathcal A}_1\cup {\mathcal A}_2}$ into a Hamilton cycle. These longer paths \jtb{found in (iii)} only use block edges.

The property \col\ will be defined differently from previous cases. We will split \col\ into three parts. Let \col1 denote the property that
\[
\mbox{$|F\cap \calN_{G}(v)|\le d_G(v)-2$ for all $v\in G$ such that $d_G(v)\ge 2$.}
\]
where $F=E(G)\setminus E(G_b)$. With a simple first moment argument we can prove the following.
\begin{claim}\lab{claim:col1}
A.a.s.\ $(G,G_b)\in\col1$.
\end{claim}

\begin{figure}
  \[
  \includegraphics[scale=0.6]{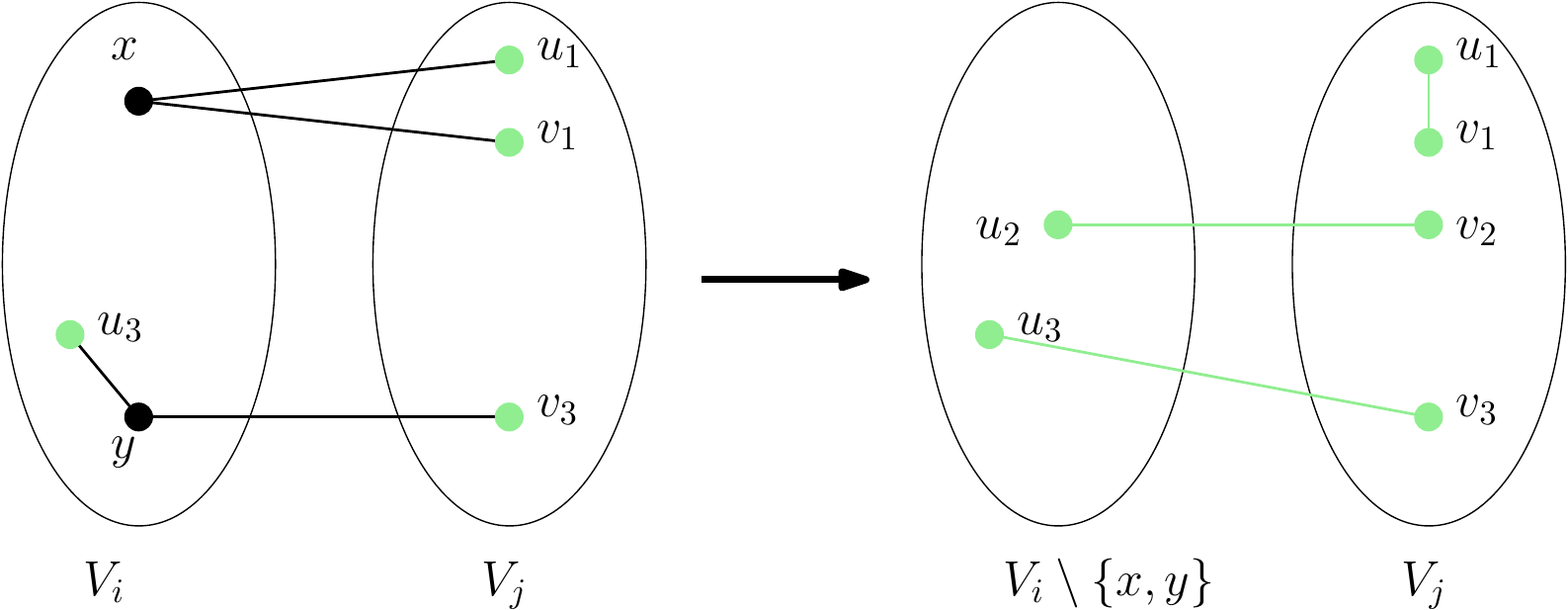}  
  \]
  \caption{Construct $M$}
  \label{f:f1}
\end{figure}

 Let $\col=\col1\cap \col2\cap \col3$ where \col2\ and \col3\ will be defined later. Note that
 $\{G\in \degree \}\cap \col1$ implies $G_b\in \degree$. If $\{G\in \degree \}\cap \col1$ holds then every problematic vertex is incident with at least 2 edges. For every $u\in {\mathcal P}$,
randomly choose 2 edges incident with $u$, colour them green. For each green edge, colour the end other than the problematic vertex green. \jt{See the left side of Figure~\ref{f:f1} for an example.}
 
 Given a path $u_1u_2\ldots u_{\ell}$, we say that we {\em supplant} the path by an edge $e=u_1u_{\ell}$ if we \jt{delete all the internal vertices on the path and their incident edges, and add edge $e$}.
Assume $\{G\in \degree \}\cap \col1$. Supplant every green 2-path in $G$ by a new green edge. Call the resulting graph $H$. Note that  $H$ is not defined if $\{G\in \degree \}\cap \col1$ fails. Note also that $H[V_i\setminus \calP]=G[V_i\setminus \calP]$ for every $1\le i\le k$. Let $E_0$ denote the set of green edges  and let $U_0$ denote the set of green vertices obtained so far in $H$.

 Next, we will choose a set of blue crossing edges and recolour them green, and colour the ends of these edges green.
 For every $1\le i<j\le k$, 
 if there are an odd number of green edges between $V_i$ and $V_j$ in $H$, then randomly choose a blue crossing edge $x$ between $V_i$ and $V_j$ in $G_b$ and recolour it green. Colour the end vertices of $x$ green.  If there is no green edge between $V_i$ and $V_j$ in $H$, then randomly choose two blue crossing edges $x,y$ between $V_i$ and $V_j$ in $G_b$ and recolour them green. Colour the end vertices of $x$ and $y$ green.

See Figure~\ref{f:f1}  for an illustration of the construction of $H$ and $E_g$.
Let $E$ denote the set of crossing edges recoloured from blue to green. Let $U$ denote the set of end vertices of edges in $E$.  Write $E=\perp$ if the above construction cannot be completed. This happens only if $|E_{G_b}(V_i,V_j)|<2$ for some $i\neq j$. However (A2) and (A4) ensure that $\mathbb{P} (|E_{G_b}(V_i,V_j)|<2)=o(1)$. The following a.a.s.\ properties are straightforward and we omit their proofs.
\begin{claim} \label{claim:E0}
A.a.s.\ the following statements hold.
\begin{enumerate}
\item[(a)] $E\neq \perp$.
\item[(b)] $E\cap E_0=\emptyset$. 
\item[(c)] $E$ induces a matching.
\item[(d)] $U\cap U_0=\emptyset$. 
\end{enumerate}
\end{claim}

 Let $E_g=E_0\cup E$ denote
all the green edges obtained so far in $H$. 
\begin{claim}
If  $H$ has a Hamiltonian cycle containing all edges in $E_g$ then $G$ contains a Hamilton cycle.
\end{claim} 
\proof Let $C$ be such a Hamilton cycle in $H$. Replacing each edge $x\in E_0\cap C$ by \jt{the} green 2-path whose supplantation yielded $x$ gives a Hamilton cycle in $G$.\qed\smallskip

We will prove that $H$ has a Hamiltonian cycle containing all edges in $E_g$ by using the following lemma. Given a set of edges $B$, let $V(B)$ denote the set of vertices spanned by $B$. Let $B_{i,j}$,  $i\leq j$, denote the set of pairs of vertices of $B$ with one end in $i$ and the other end in $j$. \jt{Given a graph $G$ and a set $E$ of edges on $V(G)$, let $G+E$ denote the graph on $V(G)$ obtained by taking the union of the edges from $G$ and $E$. }
\begin{lemma}\lab{lemma:Hcycle}
Let $B$ be a set of pairs of vertices such that
\begin{itemize}
\item the pairs in $B$ are pairwise disjoint;
\item $|B_{i,j}|>0$ for every $1\le i<j\le k$ and $\sum_{\jt{j\neq i}} |B_{i,j}|$ is even for every $1\le i\le k$;
\item $V(B)=O(\log n)$;
\item $E_g\subseteq B$;
\item no two vertices in $V(B)$ share a common neighbour in $H$;
\item no vertex in $V(B)$ is adjacent to a vertex with degree at most 2 in $H$.
\end{itemize}
Then a.a.s.\ if $\{G\in\degree\}\cap\col1$ holds then $H+B$ has a Hamilton cycle that contains all edges in $B$.
\end{lemma}
Using Claim~\ref{claim:E0} it is easy to show that $B=E_g$ satisfies all assumptions of Lemma~\ref{lemma:Hcycle}. Let ${\tt HHAM}$ denote the event that $H$ has a Hamilton cycle containing all edges in $E_g$.
If  $\{G\in\degree\}\cap \col1$ holds  then taking $B=E_g$ in Lemma~\ref{lemma:Hcycle} immediately implies {\tt HHAM}, which gives
\begin{align*}
\pr(G\in \ham)&\ge \pr(\{G\in\degree\}\cap\col1\cap {\tt HHAM} )\\
&\ge\pr(G\in\degree)-\pr(\{G\in\degree\}\cap\col1\cap \overline{{\tt HHAM}}\})-\pr(\overline{\col1})\\
&=\pr(G\in\degree)-o(1),
\end{align*}
by Claim~\ref{claim:col1} and Lemma~\ref{lemma:Hcycle}, which completes the proof of Theorem~\ref{thm:ham} in the case $q=O(1/n)$. \qed


It only remains to prove Lemma~\ref{lemma:Hcycle}. We will prove it by induction on $k$.
The following key lemma will be used to complete the inductive argument.
\begin{lemma} \label{lem:key}
Fix $1\le i\le k$. Let $A_i$ be a set of pairs of vertices  of $H_i$ such that 
\begin{itemize}
\item the pairs in $A_i$ are pairwise disjoint;
\item $V(A_i)\le \log n$;
\item no two vertices in $V(A_i)$ share a common neighbour in $H_i$;
\item no vertex in $V(A_i)$ is adjacent to a vertex with degree at most 2 in $H_i$.
\end{itemize}
Then a.a.s.\ if $\{G\in\degree\}\cap\col1$ holds then $H_i+A_i$ has a Hamilton cycle containing all of the edges in $A_i$.
\end{lemma}

We will prove  Lemma~\ref{lemma:Hcycle} in Section~\ref{sec:induction} and prove Lemma~\ref{lem:key} in Section~\ref{sec:key}. 

\subsection{Proof of  Lemma~\ref{lemma:Hcycle}}
\label{sec:induction}

We proceed with induction on $k$. The base case $k=1$ follows by Lemma~\ref{lem:key}.  Assume $k\ge 2$ and that the assertion holds for $k-1$.  

\begin{figure}
  \[
  \includegraphics[scale=0.8]{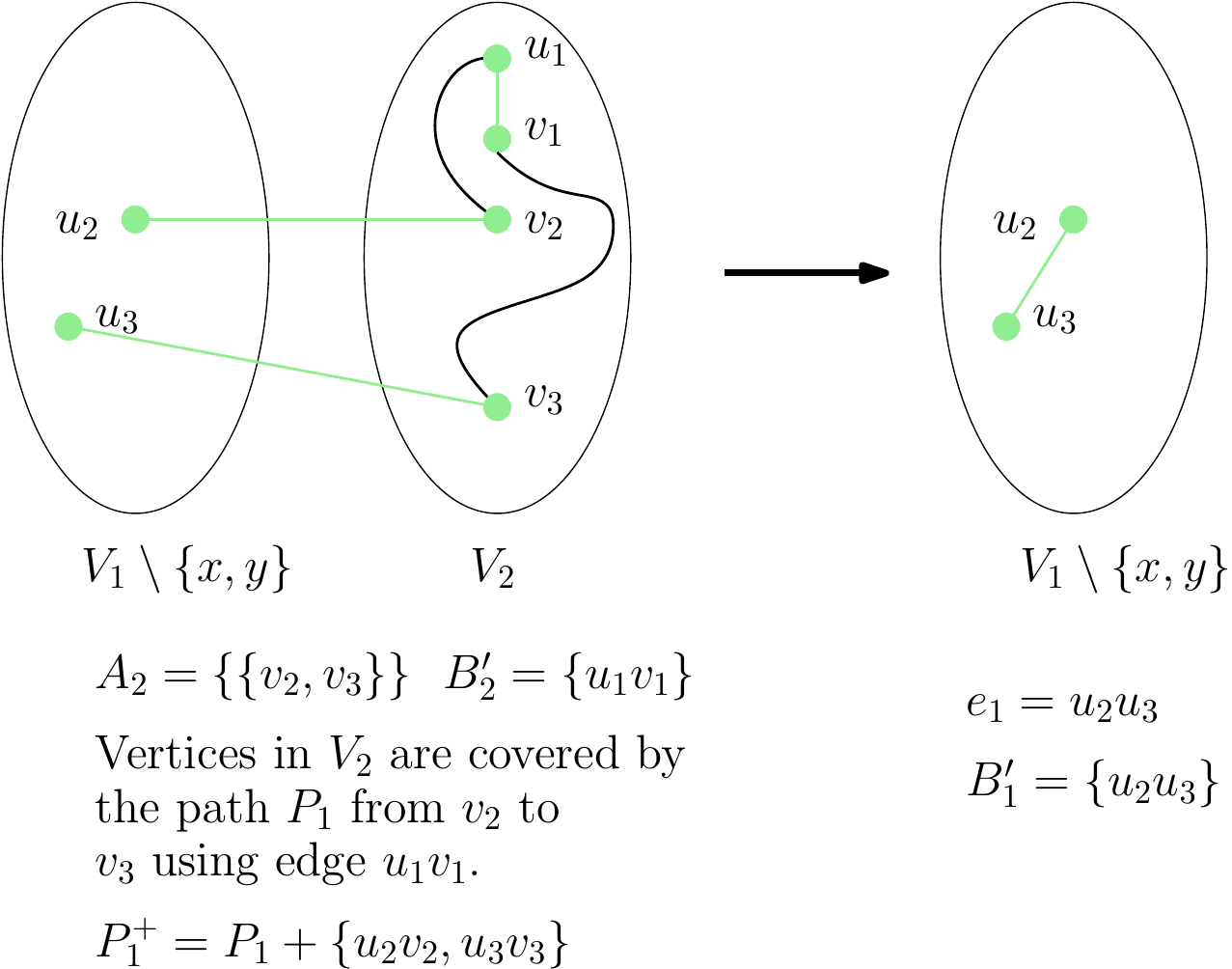}  
  \]
  \caption{Construct $P_j$, $P_j^+$ and $e_j$ with $k=2$}
  \label{f:f2}
\end{figure}

 With a slight abuse of notation we call the pairs in $B$ edges, even though they are not necessarily edges present in $H$. Let $B'_k$ denote the set of edges in $B$ with both ends in $V_k$
and let $B''_k$ denote the set of edges of $B$ with exactly one end in $V_k$. Let $V_k(B''_k)$ denote the ends of the edges in $B''_k$ that are in $V_k$. The second assumption of Lemma~\ref{lemma:Hcycle} implies that $|B''_k|$ is even. Take an arbitrary pairing $A_k'$ of the vertices in $V_k(B''_k)$ and let $A_k=B_k' \cup A_k'$. By Lemma~\ref{lem:key}, $H_k+A_k$ has a Hamilton cycle $C$ which uses all edges in $A_k$. Delete all edges in $A_k'$ from $C$. This results in a collection of vertex disjoint paths $P_1,\ldots, P_{\ell}$ such that
\begin{itemize}
\item the $\ell$ paths cover all vertices in $H_k$ and use all the edges in $B'_k$;
\item $\ell=|A_k'|=\frac{1}{2}|V_k(B''_k)|$;
\item the ends of the $\ell$ paths are the set of vertices in $V_k(B''_k)$;
\end{itemize} 
For every $P_j$ above, the ends of $P_j$ are each incident with an edge in $E_g\subseteq B''_k$. Let $P_j^+$ denote the path obtained by adding these two edges to $P_j$. Supplant $P_j^+$ by a new edge $e_j$. Now both ends of $e_j$ are in $\cup_{i\le k-1} V_i$.  See Figure~\ref{f:f2} for an example of the construction of $P_j$, $P_j^+$ and $e_j$.

Let $\widehat B=(B\cup \{e_1,\ldots, e_{\ell}\})\setminus (B'_k\cup B''_k)$. Let $\widehat H= H[\cup_{i\le k-1} V_i]$ be the subgraph of $H$ induced by $\cup_{i\le k-1} V_i$. Now $\widehat H$ has $\widehat{k}=k-1$ blocks of vertices, and $\widehat B$ is a set of pairs of vertices with both ends in $\widehat H$. Moreover, all assumptions of Lemma~\ref{lemma:Hcycle} are satisfied by $\widehat B$ and $\widehat H$ with $k$ replaced by $\widehat{k}$. By the induction hypothesis, there is a Hamilton cycle $C'$ in $\widehat H+\widehat B$ which uses all the edges in $\widehat B$. Now replacing every edge $e_j$ in $C'$ by $P_j^+$. The resulting is a Hamilton cycle in $H+B$ which uses all of the edges in $B$.

Now Lemma~\ref{lemma:Hcycle} follows by induction. \qed


\subsection{Proof of Lemma~\ref{lem:key}}
\label{sec:key}

 We introduce/recall the following useful notations.
\begin{itemize}
\item $H_i=H[V_i\setminus {\mathcal P}]=G[V_i\setminus {\mathcal P}]$;
\item $H_i^+=H_i+A_i$;
\item $H_i^{b}=G_b[V_i\setminus {\mathcal P}]$;
\item $H_i^{b+}=H_i^{b}+A_i$.
\end{itemize}

We say a Hamilton cycle in $H_i^+$ is admissible if
 it uses all edges in $A_i$.
Let $H_i^+\in\aham$ denote the property that $H_i^+$ contains an admissible Hamilton cycle.  It is sufficient to prove that 
\begin{equation}
\pr(\{G\in\degree\}\cap\col1\cap \{H_i^+\notin\aham\})=o(1). \lab{eq:prob}
\end{equation}

We will use P\'osa rotations to bound the above probability. A path in $H_i^+$ ($H_i^{b+}$) is said admissible if
 for each edge in $A_i$ it uses either that edge or none of its  vertices.
Let $L(H^+_i)$ and $L(H_i^{b+})$ denote the length of a longest admissible path in $H^+_i$ and $H_i^{b+}$ respectively. We will adapt the previous P\'osa rotation arguments to cope with admissible paths. This requires modifications of several previous definitions.  

\subsubsection{\col,\ \expanddd,\ \typical\ and \length}

First, we define \expanddd\ which is a stronger condition than \expandd. For $V\subseteq V(G)$, we say that $(G,V)$ has property $\expand^+$, if there exists an absolute constant $\eps_0>0$ such that
\[
\mbox{for every $S\subseteq V(G)$ where $|S|\le \eps_0 n$, $|N_G(S)\setminus V|\ge 2|S|\cdot \ind{n_1(G)=0}$}.
\] 
We say $G$ has property \expanddd, if the following holds: For any $F\subseteq E(G)$ such that $|F\cap {\mathcal N}_G(v)|$ is 0 if $v$ is small, and is at most $\log n/100$ if $v$ is large, and for any $V\subseteq V(G)$ such that $|V|\le \log n$, no two vertices in $V$ share a common neighbour, and no vertex in $V$ is adjacent to a vertex with degree at most 2,  \jt{we have that} $G-F$ is connected, and $(G-F,V)\in\expand^+$.
We have the following lemma whose proof is postponed until Section~\ref{sec:technical}.
\begin{lemma}\label{lem:newT1} \jt{A.a.s.\ $H_i\in\text{\expanddd}$ for every $1\le i\le k$.}
\end{lemma}
Let \col2(i) denote the event that
$L(H_i^{+})=L(H_i^{b+})$ if $L(H_i^{b+})<n-1$, and $H_i^{+}\notin \aham$ if $L(H_i^{b+})=n-1$. We may define $\col2=\cup_{i\in[k]} \col2(i)$, although in the proof of Lemma~\ref{lem:key} we only need to consider $\col2(i)$. Let \col3 denote the event that
\[
\mbox{$|F\cap \calN_{G}(v)|$ is 0 if $v$ is small in $G$, and is at most $\log n/100$ if $v$ is large in $G$,}
\]
where $F=E(G\setminus G_b)$.

We redefine \typical\ so that (T1) is replaced by 
\[
\jt{\mbox{(T1'):\quad
$H_i\in\text{\expanddd}$ for every $1\le i\le k$.}}
\]

Let $H_i^+\in \length$ denote the event that the longest admissible path  has the same number of vertices as the longest admissible cycle in $H_i^+$.

\subsubsection{Completing the proof of Lemma~\ref{lem:key}}

Because $A_i$ is a set of vertex-disjoint edges, no two edges in $A_i$ appear next to each other in the longest admissible cycle. Then, if $H_i^+$ is connected, one can always extend a longest admissible cycle to a longer admissible path, unless $H_i^+\in\aham$. Hence, $H_i^+\in\length\cap\connect$ implies that $H_i^+\in\aham$. Recall again that $H_i$ and $H_i^+$ are defined only if $\{G\in\degree\}\cap\col1$ holds. Hence
\begin{align*}
&\pr(G\in\degree\cap\col1\cap H_i^+\notin\aham)\le \pr(G\in\degree\cap\col1 \cap \{H_i^+\notin\length\cap\connect\})\\
&\quad\le \pr(G\in\degree\cap\col1 \cap H_i^+\notin\length)+\pr(G\in\degree\cap\col1 \cap H_i^+\notin\connect)\\
&\quad=\pr(G\in\degree\cap\col1 \cap H_i^+\notin\length)+o(1) \quad \mbox{(by Lemma~\ref{lem:newT1}(a)).}
\end{align*}
It is sufficient to prove that $\pr(G\in\degree\cap\col1 \cap H_i^+\notin\length)=o(1)$, which follows from the following two lemmas and the fact that a.a.s.\ $G\in\typical$.
\begin{lemma}\lab{lem:colour3}
$\pr(\col2(i)\cap \col3\mid H_i^+\in \overline{\length}\cap \{G\in \degree\cap\col1\cap \typical\})\ge \exp(-O(an/\log^2 n))$.
\end{lemma}
The proof is almost identical to the proof of Lemma~\ref{lem:colour}, with a few trivial modifications as in Lemma~\ref{lem:colour2}. We omit the details.

\begin{lemma}\lab{lem:length3}
$\pr(H_i^+\in\overline{\length}\cap \{G\in \degree\cap\col1\cap\typical\}\cap\col2(i)\cap\col3)\le \exp(-\Omega(an/\log n))$.
\end{lemma}
\proof \jt{Recall that $V(A_i)$ denotes the set of vertices spanned by the edges in $A_i$.} We have the following claim similar to Claim~\ref{claim:imply}.
\begin{claim}\label{claim:expand}
$\{G\in \degree\cap\col3\cap \typical\}$ implies that $\{H_i^{b+}\in \connect\cap \degree\}\cap\{(H_i^{b+},V(A_i))\in\expand^+\}$.
\end{claim}
Hence,
\begin{align*}
&\pr(H_i^+\in \overline{\length}\cap (G\in \degree\cap\col1\cap \typical)\cap \col2(i)\cap \col3)\\
&\quad \le \pr(H_i^+\in \overline{\length}\cap \{H_i^{b+}\in \connect\cap \degree\}\cap\{(H_i^{b+},V(A_i))\in\expand^+\}\cap \col2(i)\cap\col3)\\
&\quad \le \pr(\col2(i)\cap\col3\mid \{H_i^{b+}\in \connect\cap \degree\}\cap(H_i^{b+},V(A_i))\in\expand^+)
\end{align*}

Let $P=v_0v_1,\ldots v_{\ell}$ be a longest admissible  path in $H_i^{b+}$. A P\'osa rotation which adds edge $v_hv_{\ell}$ and deletes edge $v_hv_{h+1}$ is said to be admissible if $v_hv_{h+1}\notin A_i$.   Let $\text{End}(v_0)$ be the set of admissible paths obtained by doing admissible P\'osa rotations on $P$. We first show that 

\begin{claim}\label{cl1}
$|N_{H_i^{b+}}(\text{End}(v_0))\setminus V(A_i)|<2|\text{End}(v_0)|$.
\end{claim}
\proof The proof is similar to the standard P\'osa rotation argument. Consider any $y\in \text{End}(v_0)$ and the path $P'$  obtained via a P\'osa rotation when $y$ is added to $\text{End}(v_0)$. Assume $xy$ is an edge where $x$ in on $P$. Assume $x=v_i$ and assume $x\notin V(A_i)$. Then, either the two neightbours of $x$ on $P'$ are exactly $v_{i-1}$ and $v_{i+1}$, in which case one of them can be added to $\text{End}(v_0)$ by a P\'osa rotation; or the two neighbours of $x$ on $P'$ are not  $v_{i-1}$ and $v_{i+1}$, which implies that one of them must have been added to $\text{End}(v_0)$ before $y$. Hence, either $\{v_{i-1},v_{i+1}\}\cap\text{End}(v_0)\neq \emptyset$ or $v_i\in V(A_i)$. Our claim follows immediately.\qed \ss

By (T1') and Claim \ref{cl1}, we have that $|\text{End}(v_0)|=\Omega(n)$. Take any $x\in \text{End}(v_0)$, consider $\text{End}(x)$, the set of longest admissible paths starting from $x$ by performing admissible P\'osa rotations. Then we also have $|\text{End}(x)|=\Omega(n)$ for every $x\in\text{End}(v_0)$. If any of the edges in $E(H_i^{+}\setminus H_i^{b+})$ belongs to the set $\sE:=\{xy: x\in \text{End}(v_0),\ y\in \text{End}(x)\}$, then the event \col2(i) fails. As $|\sE|=\Omega(n^2)$, the probability that $E(H_i^{+}\setminus H_i^{b+})\cap \sE=\emptyset$ is at most $(1-\bar p)^{\Omega(n^2)}=\exp(-\Omega(an/\log n))$. Hence 
\[
\pr(\col2(i)\mid H_i^{b+}\in \connect\cap \degree\cap((H_i^{b+},V(A_i))\in\expand^+)) \le \exp(-\Omega(an/\log n)),
\]
completing the proof. \qed

\section{Proof of Lemmas~\ref{lem:expansion} and~\ref{lem:newT1}}
\lab{sec:technical}
We first state a lemma, whose proof is sketched below.
\begin{lemma}\lab{lem:properties}
Assume (A1).  A.a.s.\ $\G(\bfn,p,q)$ satisfies the following graph properties.
\begin{align*}
(C1) &\quad \text{For every $i\in[k]$, at most $n^{0.9}$ vertices in $V_i$ have degree less than $\frac{1}{2}\log n$.}\\
(C2)&\quad \text{No two vertices with degree less than 100 are within distance 5}.\\
(C3) &\quad \text{No vertices with degree less than 100 are contained in cycles of length at most 5}.\\
(C4) &\quad\text{Every set $S$ with $|S|<n/\log^2 n$ induces at most $3|S|$ edges}.\\
(C5) &\quad \text{For all $\eps>0$ there exists $\delta>0$ such that for all $S$ where $\Omega(n/\log^2 n)=|S|\le \delta n$, }\\
&\quad\text{$|E(S)|<\eps|S|\log n$.}
\end{align*}
\end{lemma}

{\em Proof of Lemma~\ref{lem:expansion}.\ }  Let $G$ be a graph satisfying properties (C1)--(C5). We also assume that the minimum degree of $G$ is at least 2, as otherwise the lemma is trivially true. In the proof of the lemma we consider various ranges for $|S|$. Colour the edges in $F$ red and let $G'=G-F$. Our assumption on $F$ implies that 
\be
\mbox{every vertex is incident with at most $\log n/100$ red edges.}\lab{F}
\ee  
Let $\eps=1/24$, and let $\delta>0$ be the constant in (C5). 

{\em Case a:} $ n/(\log n)^2\le |S| \le (\delta/3) n$. 
Let $E_1=E_G(S,\overline{S})$ and $E_2=E_{G'}(S,\overline{S})$ and let $U={\mathcal N}_{G'}(S)$. Suppose that $|U|<2|S|$. Then, $|S\cup U| <3|S|\le \delta n$. By (C5), $S\cup U$ induces at most $\eps|S\cup U|\log n\le 3\eps |S|\log n$ edges in $G$. This implies that $|E_2|\le 3\eps|S|\log n$. By (C1), the total degree of vertices in $S$ is at least $(|S|-n^{0.9})\cdot (1/2)\log n\ge |S|\log n/3$. On the other hand, by (C5), $S$ induces at most $\eps |S|\log n$ edges. Thus, $|E_1|\ge |S|\log n/3- 2\eps |S|\log n = |S|\log n/4$. Consequently,  $|F\cap E_G(S,\overline{S})|=|E_1|-|E_2|\ge (5/24)|S|\log n$, contradicting condition~\eqn{F}.  

{\em Case b} $|S|\le n/\log^2 n$.
A vertex in $G$ is called {\em extremely small} if its degree is less than 100. 
 Let $\calX$ denote the set of extremely small vertices in $S$, and $\calY$ denote $S\setminus \calX$. 

{\em Case b1: $|\calY|=0$.} Then $F$ is not incident with any vertex in $S$. By (C2), $S$ must induce an independent set. By our assumption, all vertices in $S$ has degree at least 2. By (C2), $\calN_{G}(a)\cap \calN_{G}(b)=\emptyset$ for every distinct $a,b\in S$. It follows immediately then that $|\calN_{G'}(S)|=|\calN_G(S)|\ge 2|S|$.

{\em Case b2: $|\calY|\ge 1$.} Now $F$ is not incident with any vertex in $\calX$. Let $\calZ_1={\mathcal N}_{G'}(\calX)\setminus \calY={\mathcal N}_{G}(\calX)\setminus \calY$ be the set of neighbours of $\calX$ that are not in $\calY$, $\calZ_2=N_{G'}(S)\setminus\calZ_1$ be the neighbours of $S$ in $G'$ that are not in $\calZ_1$. Then, $|\calN_{G'}(S)|=|\calZ_1|+|\calZ_2|$. Let $\calY_1=\calN_{G'}(\calX)\cap \calY$ be the set of neighbours of $\calX$ in $\calY$. By (C2) and our assumption that the minimum degree of $G$ is at least 2, 
\begin{align}
|\calY_1|&=|E_{G'}(\calX,\calY)|=|E_{G}(\calX,\calY)|\lab{c1} \\
|\calN_{G'}(\calX)|&=|\calZ_1|+|\calY_1|\ge 2|\calX|. \lab{c2}
\end{align}
We prove next that every vertex in $\calY$ can be incident to at most one vertex in $\calZ_1$ in $G$ (and $G'$). Assume $a\in \calY$ is adjacent to two vertices $b$ and $c$ in $\calZ_1$. If $b$ and $c$ have a common neighbour $z\in \calX$, then $abzc$ forms a 4-cycle in $G$, violating (C3). Assume $b$ and $c$ each adjacent to $b'\in\calX$ and $c'\in\calX$ respectively. Then $b'bacc'$ is a 4-path in $G$ connecting two light vertices, violating (C2). Hence, $|\calN_{G'}(z)\cap \calZ_1|\le1$ for every $z\in \calY$. Consequently, 
\be
|E_{G'}(\calY,\calZ_1)|\le |\calY|.\lab{c3}
\ee

Assume to the contrary that $|\calZ_1|+|\calZ_2|=|\calN_{G'}(S)|<2|S|=2(|\calX|+|\calY|)$. Then, by~\eqn{c2} we have
\be
|\calZ_2|<2|\calY|+|\calY_1|.\lab{c4}
\ee
Every vertex in $\calY$ has degree at least 100 in $G$. Also, a vertex is incident with red edges only if its degree is at least $\log n/10$ and then at most $\log n/100$ red edges. It follows that every vertex in $\calY$ has degree at least 100 in $G'$ as well. Thus,
\be
|E_{G'}(\calX,\calY)|+2|E_{G'}(\calY)|+|E_{G'}(\calY,\calZ_1)|+|E_{G'}(\calY,\calZ_2)|\ge 100|\calY|. \lab{c5}
\ee
By (C4), $|E_{G'}(\calY)|\le 3|\calY|$, $|E_{G'}(\calY,\calZ_2)|\le |E_{G'}(\calY\cup \calZ_2)|\le 3(|\calY|+|\calZ_2|)$. We have shown that $|E_{G'}(\calX,\calY)|\le |\calY_1|$ and $|E_{G'}(\calY,\calZ_1)|\le |\calY|$. Hence, by~\eqn{c4}, the left hand side of~\eqn{c5} is at most $16|\calY|+4|\calY_1|\le 20|\calY|$, whereas the right hand side is $100|\calY|$, contradiction. This confirms that $|N_{G'}(S)|\ge 2|S|$.\qed\ss

\no {\em Proof of Lemma~\ref{lem:properties}. } The proof is standard and straightforward. We give a sketch only and omit the somewhat tedious calculations. For (C1), following the same argument as in Lemma~\ref{lem:mindegree}, we can show that the expected number of vertices with degree less than $\frac{1}{2}\log n$ is $o(n^{0.9})$.

By (A1), there is $i\in[k]$ such that $pn_i+(n-n_i)q -\log n_i=\log\log n+O(1)$. Together with (A4), this implies that $p,q=O(\log n/n)$. Using this, the expected number of $S$ with $|S|=s$ which induce more than $3s$ edges is at most
\[
\binom{n}{s} \binom{s^2}{3s} (C\log n/n)^{3s}, \quad \mbox{for some constant $C>0$}.
\]
It is straightforward to see that summing the above over $s<n/\log^2n$ yields $o(1)$. This proves (C4). The proof for (C5) is almost the same, and we can bound the expected number of sets $S$ inducing more than $\eps|S|\log n$ by $o(1)$ by choosing a sufficiently small $\delta$.

For (C2), we bound the expected number of $\ell$-paths, $\ell\le 5$, where the ends are vertices of degree less than 100. There are at most $n^{\ell+1}$ ways to choose the $\ell+1$ vertices.  Using the probability bounds as in Lemma~\ref{lem:mindegree}, the probability of both of the chosen end-vertices having degree less than 100 is at most $\log^{200}n/n^2$. The probability that the chosen $\ell+1$ vertices form a path is bounded by $(C\log n/n)^{\ell}$. Multiplying all together we have that the expected number of such paths is at most $n^{\ell+1}\cdot (\log^{200}n/n^2)\cdot  (C\log n/n)^{\ell}=o(1)$. This proves (C2). The proof of (C3) is similar.
\qed\ss

\no{\em Proof of Lemma~\ref{lem:newT1}.\ } \remove{
There are three types of vertices in $V(A_i)$: (a) they are the neighbours of some vertex in $V_i$ whose degree in $G[V_i]$ equals 1; (b) they are the neighbours of some vertex $V_j$ (for some $j\neq i$) whose degree in $G[V_j]$  is at most 1; (c) they are in $U$, i.e.\ they are ends of some crossing edge recoloured from blue to green. Vertices of the first two types are neighbours of vertices with degree at most 1, and vertices of the last type are random vertices in $V_i$. It is standard argument to show $|V(A_i)|\le \log n$, and no two vertices in $V(A_i)$ share a common neighbour, and no vertex in $V(A_i)$ is adjacent to a vertex of degree at most 2.  That confirms part (b).
}
We only need to consider $\G(n_i,p)$ where $pn_i\ge \log n+\log\log n +O(1)$. The proof is almost the same as the proof of Lemma~\ref{lem:expansion} with only small modifications which we point out below. Again we assume that $H_i$ is a graph satisfying (C1)--(C5) (for (C1), we only need to consider a fixed $i$). Let $V$ be an arbitrary set of vertices in $H_i$ such that $|V|\le \log n$, no two vertices in $V$ share  a common neighbour, and no vertex in $V$ is adjacent to a vertex of degree at most 2. Colour the vertices in $V$ red. Let $F$ be an arbitrary set of edges such that $|F\cap \calN_{H_i}(v)|$ is 0 if $v$ is small, and is at most $\log n/100$ if $v$ is large. Let $H_i'=H_i-F$.

Let $\eps=1/24$, and let $\delta>0$ be the constant in (C5). 
For $S$ where $ n/(\log n)^2\le |S| \le (\delta/4) n$, let $E_1=E_{H_i}(S,\overline{S})$ and let $E_2=E_{H_i'}(S,\overline{S})$. Let $U={\mathcal N}_{H_i'}(S)$. Assume $|U\setminus V|<2|S|$. Then, $|S\cup U| <3|S|+|V|\le 4|S|\le \delta n$. Now with the same proof as in Lemma~\ref{lem:expansion}, we can lead to a contradiction with condition~\eqn{F}. Thus, we must have $|U\setminus V|\ge 2|S|$ in this case.
 
For $S$ where $|S|\le n/\log^2 n$, 
again,
call a vertex in $H_i$ {\em extremely small} if its degree is less than 100. 
 Let $\calX$ denote the set of extremely small vertices in $S$, and $\calY$ denote $S\setminus \calX$. 

{\em Case 1: $|\calY|=0$.} Then $F$ is not incident with any vertex in $S$. By (C2), $S$ must induce an independent set. Recall that $H_i=G_b[V_i\setminus \calP]$. It follows immediately that every vertex in $H_i$ has degree at least 2, as all problematic vertices are in $\calP$. Moreover, by the assumptions on $V$, the vertices with degree 2 are not adjacent to any vertex in $V$, and every vertex in $H_i$ can be adjacent to at most one vertex in $V$. Hence, every vertex has at least 2 non-red neighbours. By (C2), $\calN_{H_i}(a)\cap \calN_{H_i}(b)=\emptyset$ for every distinct $a,b\in S$. It follows immediately then that $|\calN_{H_i'}(S)\setminus V|=|\calN_{H_i}(S)\setminus V|\ge 2|S|$.

{\em Case 2: $|\calY|\ge 1$.} Now $F$ is not incident with any vertex in $\calX$. Let $\calZ_1={\mathcal N}_{G'}(\calX)\setminus \calY={\mathcal N}_{G}(\calX)\setminus \calY$ be the set of neighbours of $\calX$ that are not in $\calY$, $\calZ_2=N_{G'}(S)\setminus\calZ_1$ be the neighbours of $S$ in $G'$ that are not counted in $\calZ_1$. Let $\calZ_i'=\calZ_i\setminus V$ for $i\in\{1,2\}$.  Then, $|\calN_{H_i'}(S)\setminus V|=|\calZ'_1|+|\calZ'_2|$. Let $\calY_1=\calN_{G'}(\calX)\cap \calY$ be the set of neighbours of $\calX$ in $\calY$. With the same argument for~\eqn{c1}--\eqn{c3}, together with the fact that every vertex in $\calX$ has at least 2 non-red neighbours, we have
\begin{align}
|\calY_1|&=|E_{H_i'}(\calX,\calY)|=|E_{H_i}(\calX,\calY)|\lab{c1'} \\
|\calN_{H_i}(\calX)\setminus V|&=|\calZ_1'|+|\calY_1|\ge 2|\calX| \lab{c2'}\\
|E_{H_i'}(\calY,\calZ_1)|&\le |\calY|.\lab{c3'}
\end{align}
Assume to the contrary that $|\calZ_1'|+|\calZ_2'|=|\calN_{H_i'}(S)\setminus V|<2|S|=2(|\calX|+|\calY|)$. Then, by~\eqn{c2'} we have
\be
|\calZ_2'|<2|\calY|+|\calY_1|.\lab{c4'}
\ee
Every vertex is adjacent to at most one vertex in $V$. It follows immediately that
\be
|\calZ_2|\le|\calZ_2'|+|\calY|. \lab{Z2}
\ee  
Every vertex in $\calY$ has degree at least 100 in $H_i$. Also, a vertex is incident with red edges only if its degree is at least $\log n/10$ and then at most $\log n/100$ red edges. It follows that every vertex in $\calY$ has degree at least 100 in $G'$ as well. Thus,
\be
|E_{H_i'}(\calX,\calY)|+2|E_{H_i'}(\calY)|+|E_{H_i'}(\calY,\calZ_1)|+|E_{H_i'}(\calY,\calZ_2)|\ge 99|\calY|. \lab{c5'}
\ee
By (C4), $|E_{H_i'}(\calY)|\le 3|\calY|$, $|E_{H_i'}(\calY,\calZ_2)|\le |E_{H_i'}(\calY\cup \calZ_2)|\le 3(|\calY|+|\calZ_2|)\le 3(4|\calY|+|\calY_1|)$ by~\eqn{c4'} and~\eqn{Z2}. We have shown that $|E_{H_i'}(\calX,\calY)|=|\calY_1|$ and $|E_{H_i'}(\calY,\calZ_1)|\le |\calY|$. Hence,  the left hand side of~\eqn{c5} is at most $19|\calY|+4|\calY_1|\le 23|\calY|$, whereas the right hand side is $99|\calY|$, contradicting with $|\calY|\ge 1$. This shows that $|N_{H_i'}(S)\setminus V|\ge 2|S|$.\qed\ss
\section{Conclusion}
We have analysed the Hamiltonicity of a particular stochastic block model and given tight estimates for the threshold. The most natural extension of this work will be to the case where $P$ is an arbitrary symmmetric stochastic matrix. This will be the subject of further research.

\bibliographystyle{plain}

\end{document}